\documentclass[12pt]{amsart}

\newcommand{\la}{\leftarrow}
\newcommand{\stakr}{*}
\newcommand{\AC}{Arnold Conjecture}

\usepackage{mathabx}
\input bookman.sty
\boldmath

\usepackage{helvet}

\usepackage{graphics}
\usepackage{amssymb}
\usepackage{amsxtra}
\usepackage{amsmath}
\usepackage{mathrsfs}

\usepackage[arrow, matrix]{xy}
\xyoption{frame}

\theoremstyle{plain}

\newtheorem{theo}{Theorem}[section]
\newtheorem{lemm}[theo]{Lemma}
\newtheorem{prop}[theo]{Proposition}
\newtheorem{coro}[theo]{Corollary}

\theoremstyle{definition}

\newtheorem{defi}[theo]{Definition}
\newtheorem{rema}[theo]{Remark}

\newfont{\rmm}{cmr10 scaled 1000}

\newfont{\itt}{cmsl10 scaled 1000}

\newfont{\rM}{cmr10 scaled 1700}

\newcounter{lemma}[section]

\newcounter{tempcounter}

\newcommand{\lb}{\label}

\newcommand{\rrf}[1]{(\ref{#1})}

\renewcommand{\a}{\alpha}
\renewcommand{\b}{\beta}
\newcommand{\g}{\gamma}
\renewcommand{\d}{\delta}

\newcommand{\ve}{\varepsilon}

\renewcommand{\l}{\lambda}

\renewcommand{\r}{\rho}

\newcommand{\s}{\sigma}

\renewcommand{\o}{\omega}

\newcommand{\G}{\Gamma}

\renewcommand{\L}{\Lambda}

\newcommand{\EE}{{\mathcal E}}
\newcommand{\FF}{{\mathcal F}}

\newcommand{\LL}{{\mathcal L}}
\newcommand{\MM}{{\mathcal M}}

\newcommand{\PP}{{\mathcal P}}

\newcommand{\RR}{{\mathcal R}}

\newcommand{\ff}{{\mathbb{F}}}

\newcommand{\nn}{{\mathbb{N}}}

\newcommand{\qq}{{\mathbb{Q}}}
\newcommand{\rr}{{\mathbb{R}}}

\newcommand{\zz}{{\mathbb{Z}}}

\newcommand{\RRR}{{\mathbf{R}}}

\newcommand{\CCCC}{{\mathscr{C}} }

\newcommand{\Ker}{\text{\rm Ker }}

\renewcommand{\Im}{\text{\rm Im }}

\newcommand{\Id}{\text{\rm Id}}

\newcommand{\GL}{\text{\rm GL}}

\newcommand{\bere}{\begin{rema}}
\newcommand{\bede}{\begin{defi}}

\renewcommand{\beth}{\begin{theo}}
\newcommand{\bele}{\begin{lemm}}
\newcommand{\bepr}{\begin{prop}}
\newcommand{\beeq}{\begin{equation}}
\newcommand{\bega}{\begin{gather}}
\newcommand{\begaa}{\begin{gather*}}
\newcommand{\been}{\begin{enumerate}}

\newcommand{\bedee}{\begin{defii}}
\newcommand{\bethh}{\begin{theoo}}
\newcommand{\belee}{\begin{lemmm}}
\newcommand{\beprr}{\begin{propp}}

\newcommand{\beco}{\begin{coro}}

\newcommand{\beal}{\begin{aligned}}

\newcommand{\enre}{\end{rema}}

\newcommand{\enco}{\end{coro}}
\newcommand{\enpr}{\end{prop}}
\newcommand{\enth}{\end{theo}}
\newcommand{\enle}{\end{lemm}}
\newcommand{\enen}{\end{enumerate}}
\newcommand{\enga}{\end{gather}}
\newcommand{\engaa}{\end{gather*}}
\newcommand{\eneq}{\end{equation}}
\newcommand{\enal}{\end{aligned}}

\newcommand{\bq}{\begin{equation}}
\newcommand{\bqq}{\begin{equation*}}

\renewcommand{\leq}{\leqslant}
\renewcommand{\geq}{\geqslant}

\newcommand{\wi}{\widetilde}

\newcommand{\wh}{\widehat}

\newcommand{\sbs}{\subset}

\newcommand{\tens}[1]{\underset{#1}{\otimes}}

\newcommand{\wrt}{with respect to}
\newcommand{\ho}{homomorphism}

\newcommand{\ma}{manifold}

\newcommand{
\hos}{~homomorphisms}

\newcommand{\fbfg}{ free based finitely generated }

\newcommand{\Prf}{{\it Proof.\quad}}

\newcommand{\smo}{C^{\infty}}

\newcommand{\chart}{\Phi_p:U_p\to B^n(0,r_p)}
\newcommand{\atlas}{\{\Phi_p:U_p\to B^n(0,r_p)\}_{p\in S(f)}}

\newcommand{\pr}{\partial}

\newcommand{\qs}{\hfill\square}

\newcommand{\pa}{\vskip0.1in}

\newcommand{\arrh}[3]
{
\xymatrix{
{#1} \ar[r]^<<<<{#2}  &{#3}
}
}

\newcommand{\arlh}[3]
{
\xymatrix{
{#1} & \ar[l]_<<<<{#2}  {#3}
}
}

\newcommand{\arrr}[1]
{\arrh {}{#1}{}}

\newcommand{\arrl}[1]
{\arlh {}{#1}{}}

\newcommand{\arl}
{\arrl {}}

\newcommand{\arrto}
{\xymatrix{{} \ar@{|-{>}}[r]  & {} } }

\newcommand{\arrinto}
{\xymatrix{{} \ar@{^{(}->}[r]  & {} } }

\newcommand{\ffg}{free  finitely generated}
\newcommand{\cco}{ chain complex}

\begin{document}

\title
[On the fixed points of a Hamiltonian diffeomorphism]
{On the fixed points of a Hamiltonian 
diffeomorphism in presence of fundamental group}
\author{Kaoru Ono and  Andrei Pajitnov}
\address{Research Institute for Mathematical Sciences, Kyoto University, Kyoto 606-8502, Japan}
\email{ono@kurims.kyoto-u.ac.jp}
\address{Laboratoire Math\'ematiques Jean Leray 
UMR 6629,
Universit\'e de Nantes,
Facult\'e des Sciences,
2, rue de la Houssini\`ere,
44072, Nantes, Cedex}                    
\email{andrei.pajitnov@univ-nantes.fr}

\thanks{Kaoru Ono is supported by JSPS Grant-in-Aid for Scientific Research 
\#  21244002 and \# 26247006} 
\begin{abstract}
Let $M$ be a weakly monotone symplectic manifold, and $H$ be 
a time-dependent Hamiltonian; we assume that the periodic
orbits of the corresponding time-dependent 
Hamiltonian vector field 
are non-degenerate. We construct a refined version of 
the Floer 
chain complex associated to these data and any regular covering 
of $M$, and derive from it new lower bounds for the 
number of periodic orbits.

Using these invariants we prove in particular that if $\pi_1(M)$ 
is finite and solvable or simple,  then the number of 
periodic orbits is not less  than the minimal number of generators
of $\pi_1(M)$.  
For a general closed symplectic manifold with infinite fundamental group, we show the existence 
of $1$-periodic orbit of Conley-Zehnder index $1-n$ for any non-degenerate $1$-periodic Hamiltonian system.   

\end{abstract}
\maketitle

\section{Introduction}
\label{s:intro}

Let $M^{2n}$ be a closed symplectic \ma,
denote by $\o$ its symplectic form.
Let $H:\rr/\zz\times M\to\rr$ be a $\smo$ function
{\it (the Hamiltonian)}.
We will write $H_t(x)$ instead of $H(t,x)$.
One associates to $H$ a time-dependent 
Hamiltonian vector
field $X_{H_t}$ on $M$ by the formula
$$
\o(X_{H_t},\cdot)=dH_t \ \ {\rm for\  every } \ \ t.
$$

Assume that every periodic orbit of $\{ X_{H_t } \}$  is non-degenerate.  
Then the set $\PP(H)$ of all periodic orbits is finite.
Let  $p(H)$  be the cardinality of this set.
Denote by $\MM(M)$ the {\it Morse number} of $M$,  
that is, the minimal possible number of 
critical points of a Morse function on $M$.

The celebrated Arnold conjecture 
(see 
\cite{Arnold}, Appendix 9, and \cite{ArnoldProbs},
p.284) says that
\begin{equation}\label{f:gac}
p(H) \geq \MM(M).
 \end{equation}
\noindent
It was proved by V. I. Arnold himself in the case 
when $H$ is ``sufficiently small''  function.
The Arnold conjecture implies in particular 
certain homological lower bounds for
$p(H)$. Namely, 
let us denote by $b_i(M)$ the rank of $H_i(M)$,
and by $q_i(M)$ the torsion number of $H_i(M)$
(that is, the minimal possible number of generators
of the abelian group $H_i(M)$).
Then the conjecture  \rrf{f:gac} implies
the following:
\begin{equation}\label{f:zhac}
p(H) \geq \sum_i \Big(b_i(M)+q_i(M)+q_{i-1}(M)\Big).
 \end{equation}
 The inequality \rrf{f:gac} implies also the following:
 \begin{equation}\label{f:fhac}
p(H) \geq \sum_i b^\ff_i(M),
 \end{equation}
 where $\ff$ is any field and  we denote by $b^\ff_i(M)$
 the dimension of $H_i(M,\ff)$ over $\ff$.
 
 A. Floer \cite{Floer}  constructed
 a chain complex associated with a non-degenerate $1$-periodic Hamiltonian $\{ H_t \}$;
 applying this construction he proved the
 homological version \rrf{f:fhac} of the \AC~ for
 any field $\ff$  in the case of monotone symplectic manifolds.  
 (Since the degree in Floer homology is $\zz/2N \zz$, torsions in ordinary homology appearing in different degrees 
 but congruent modulo $2N$ with relatively prime orders contribute to the Floer homology in the same degree.  
 This is the reason why \rrf{f:zhac} does not follow from Floer homology with integer coefficients.
 Here $N$ is the minimal Chern number of $(M,\omega)$.)  
 
 The construction of the Floer chain complex 
 was generalized to  wider classes of symplectic manifolds, i.e., weakly monotone symplectic manifolds, 
 in \cite{Hofer-Salamon}, \cite{Ono} .  
 Taking the results on orientation \cite{Floer} section 2e, \cite{FO99} section 21 into account,  the  conjecture \rrf{f:fhac}
 is verified in the case of weakly monotone symplectic manifolds  (\cite{Floer} for monotone case, \cite{Hofer-Salamon} for $N=0$ or $N \geq n$,  \cite{Ono} 
 for weakly monotone case).   
 If the minimal Chern number is zero, i.e., spherically Calabi-Yau, the inequality \rrf{f:zhac} holds.   
The construction over  $\qq$ was further generalized to all closed symplectic manifolds in \cite{FO99}, \cite{LT}, hence 
the inequality \rrf{f:fhac}  with $\ff$ of characteristic $0$, e.g., $\ff=\qq$ follows. 
 
 These results confirm the homological versions
 of the Arnold Conjecture, i.e., \rrf{f:zhac} holds when the minimal Chern number of the closed symplectic manifold is zero, 
 \rrf{f:fhac} with any field $\ff$ holds in the case of weakly monotone closed symplectic manifolds and 
 \rrf{f:fhac} with a field of characteristic zero holds for general closed symplectic manifolds.
 As for the initial conjecture
 \rrf{f:gac} it is still unproved in general case.
 For a simply connected \ma~
 $M^{2n}$ with $n\geq 3$ the statement  
 \rrf{f:gac} is equivalent 
 to \rrf{f:zhac} in view of S. Smale's theorem
 \cite{Smale}.
 However in the non-simply-connected case 
 the number $\MM(M)$ can be strictly 
 greater than the right hand side
 of \rrf{f:zhac}. A first step to the proof of 
 the geometric Arnold Conjecture  \rrf{f:gac} would 
 be to prove a weaker inequality involving only
 the invariants of the fundamental group.
 For a group $G$ let
 \begin{equation}
 \RR = \{ 1\la G\la F_1\la F_2 \}
 \end{equation}
be a presentation of $G$, where 
$F_1$ and $F_2$ are free groups of ranks $d(\RR)$ and 
$r(\RR)$. Denote by $d(G)$ the 
minimum of numbers $d(\RR)$ for all presentations $\RR$,
and by $D(G)$ the minimum of 
numbers $d(\RR)+r(\RR)$  for all presentations $\RR$.
On the occasion of Arnoldfest in Toronto 1997, 
V. I.  Arnold asked the first author 
whether the development in Floer theory at that 
time\footnote{\cite{FO99}, \cite{LT} 
appeared as preprints in the previous year.}  
settled the original form of his conjecture, i.e., \rrf{f:gac}, and, 
in particular, whether one can show the following weaker assertion, 
which does not follow from homological version of the conjecture:
   \begin{equation}\label{f:pi1ac}
  p(H) \geq D(\pi_1(M)).
 \end{equation}
 
 A weaker form of this conjecture is the following:
  \begin{equation}\label{f:pi1ac2}
  p(H) \geq d(\pi_1(M)).
 \end{equation}
  
  Since then some progress has been made 
  in this direction, although the 
  conjecture is far from being solved.
  M. Damian 
  \cite{Damian2002} considers 
  similar questions in the framework 
  of  the Hamiltonian isotopies of the 
  cotangent bundle of a compact manifold $M$.
 In a recent preprint 
  \cite{barraud} J.-F. Barraud suggested a construction
  of  a Floer fundamental group, and proved in particular,
  that $p_{1-n}(H)\geq 1$ if $\pi_1(M)$ is non-trivial
  and $M$ is spherical Calabi-Yau or monotone.
  (Here $p_j(H)$ stands for the number of 
  periodic orbits of Conley-Zehnder index $j$.)
  
 In the present paper we use 
 the Floer chain complex associated with $H$ and a regular covering $\wi M\to M$
 of the underlying \ma,
 and deduce from it new lower bounds for $p(H)$
 in terms of certain invariants $\mu_i(\wi M)$,
 which depend on 
 the homotopy
 type of $M$, the minimal Chern number $N$
 of $M$, and the chosen covering
 (see Definition \ref{d:mui-m}). 
 These invariants  are
 similar to V. V. Sharko's 
 invariants of chain complexes
 \cite{Sharko}.
 The numbers $\mu_i$ are indexed by $\nn$ in the
 case of spherical Calabi-Yau manifolds
 and by $\zz/2N\zz$ in case when 
 the minimal Chern number of $M$ equals $N$.
 We have
 $$
 p(H)
 \geq \sum_i\mu_i(\wi M).$$
 
 Using these invariants 
 we obtain partial results in the direction of 
 the conjecture \rrf{f:pi1ac2}. For a group
 $G$ denote by $\d(G)$ the minimal number of generators
 of the augmentation ideal of
 $G$ as a $\zz[G]$-module.
 In the  case when $M$ is weakly monotone 
 and $\pi_1(M)$ is a finite group 
 we prove that
 $$
 p(H) \geq \d(\pi_1(M)).
 $$
  In particular we confirm the conjecture 
 \rrf{f:pi1ac2}
 for weakly monotone manifolds whose fundamental groups are
 finite simple or solvable (Theorem \ref{t:monoton}).  
 
 We also show the existence of $1$-periodic orbits of Conley-Zehnder index $1-n$ for any non-degenerate 
 $1$-periodic Hamiltonian system on any closed symplectic manifold  with infinite fundamental group (Theorem \ref{t:inf-gen-case}).
 
\section{Floer complex on the covering space}
\label{s:Floeroncovering}
\subsection{Novikov rings: definitions}
\label{su:def-nov}
In this  preliminary subsection we gathered 
 the definitions of several versions of the Novikov
rings with which we will be working 
in the paper.

Let $T$ be a finitely generated free abelian group,
and $\xi:T\to\rr$ be a homomorphism. 
Let $R$ be a ring (commutative or not). 
Recall that the group ring 
$R[T]$ is the set of all 
finite linear combinations
$$
l=\sum_{i=0}^N a_i g_i, \ \ \ \ \ 
{\rm with ~ ~ } a_i\in R, \ \ g_i\in T
$$
with  a natural ring structure 
(determined by the  requirement that the elements of $R$ commute
with the elements of $T$).

We denote by $R((T))$ the set of all 
formal linear combinations (infinite in general)
$$
\l=\sum_{i=0}^\infty a_i g_i, \ \ \ \ \ 
{\rm with ~ ~ } a_i\in R, \ \ g_i\in T
$$
such that $\xi(g_i)\to -\infty $ with $i\to \infty$.
Thus the series $\l$ can be infinite, but for every $C$
the number of terms of $\l$ with $\xi(g_i)\geq C$ is finite.
Usually the \ho~ $\xi$ is clear from the context, so we omit it from
the notation.
The usual definition of the product
of power series endows the abelian group $R((T))$ with
the natural ring structure
(we require that the elements of $R$ commute
with the elements of $T$). This ring is called { \it the Novikov
completion of the  ring $R[T]$}. 
In this paper we will work with the case 
$R=\zz[G]$, where $G$ is a group.

The augmentation \ho~ $\ve : \zz[G]\to\zz$ has a natural extension to a ring
\ho~ $R((T))\to\zz((T))$, which will be denoted by the same
symbol $\ve$:
$$
\ve (\sum_i a_ig_i)= \sum \ve(a_i) g_i \ \ \ \ \
{\rm with ~ ~ } a_i\in R, \ \ g_i\in T.
$$
Thus the ring $\zz((T))$ acquires a natural structure of $R((T))$-module.
\pa
\bere\label{r:fin-g}
If the group $G$ is finite, then the ring $\zz[G]((T))$
coincides with the group ring $\zz((T))[G]$ of the group 
$G$ with coefficients in $\zz((T))$.
\enre
In the case when 
$\xi:T\to\rr$ is a monomorphism,
we will use abbreviated notation. 
The group ring 
$\zz[T]$ will be denoted by $\L$, and its Novikov completion
with respect to a monomorphism $\xi$ will be denoted by $\wh\L$.
For a field $\ff$ we denote by $\FF$ the 
Novikov completion of the group ring $\ff[T]$
with respect to $\xi$. The ring $\wh\L$
is a principal ideal domain (PID), and $\FF$ is a field. 
We will denote the ring $\zz[G]((T))$ by $\LL$,
The ring $\ff[G]((T))$ will be denoted by $\LL_\ff$.
These rings will appear frequently in Sections 3 and 4.

The Novikov rings appear in Hamiltonian dynamics in the following context
(see Subsection \ref{su:review}).
Let $M$ be a closed symplectic manifold. 
The de Rham cohomology class of the symplectic form
determines a \ho~
$
[\o]:\pi_2(M)\to\rr.
$
Consider the group
$$
\Gamma 
=
\pi_2(M)/(\ker [\omega] \cap \ker c_1(M)),
$$
where $c_1(M)$ is the Chern class of the almost complex structure associated  
to $\o$.
The Novikov completion $\zz((\G))$ will be denoted 
by $\Lambda^{\mathbb Z}_{(M,\omega)}$.
The Novikov ring $\zz[G]((\G))$ will be denoted in this context by
$\Lambda^{\mathbb Z[G]}_{(M,\omega)}$.

The restriction of the \ho~ $\o$ to a smaller group
$$\Gamma_0=
\ker c_1(M)/(\ker [\omega] \cap \ker c_1(M))
$$
is a monomorphism, so the corresponding Novikov completion $\zz((\G_0))$
is a PID; it will be denoted by 
$\Lambda^{(0){\mathbb Z}}_{(M,\omega)}=\zz((\G_0))$.
The Novikov ring $\zz[G]((\G_0))$ will be denoted in this context by
$\Lambda^{(0){\mathbb Z[G]}}_{(M,\omega)}$.

\subsection{Review on Hamiltonian Floer complex}${}$
\label{su:review}

In this subsection, we recall the construction of Hamiltonian Floer complex with 
integer coefficients\footnote{There is an approach 
to construct Hamiltonian Floer complex with integer coefficients for non-degenerate 
periodic Hamiltonian systems on 
arbitrary closed symplectic manifold \cite{FO2000}.  
Since the details has not been carried out, we restrict ourselves to the class of weakly monotone 
symplectic manifolds.} following \cite{Floer},
\cite{Hofer-Salamon}, \cite{Ono}.  
Here we use homological version.  
Let $(M, \omega)$ be a closed symplectic manifold
of dimension $2n$.  
The minimal Chern number $N=N(M,\omega)$ of $(M, \omega)$ 
is a non-negative integer such that 
$\{ \langle c_1(M), A \rangle \vert A \in \pi_2(M) \} = N {\mathbb Z}$.  
We call $(M, \omega)$ weakly monotone (semi-positive) if 
$\langle [\omega], A \rangle \leq 0$ holds for any $A \in \pi_2(M)$ 
with $3-n \leq \langle c_1(M), A \rangle < 0$.  
This class of symplectic manifolds, in particular,  
contains the following.  

\been
\item (monotone case) We call $(M,\omega)$ a monotone symplectic manifold,
if there exists a positive real number $\lambda$ such that the following 
equality holds 
$$\langle c_1(M), A \rangle = \lambda \langle [\omega], A \rangle$$ 
for any $A \in \pi_2(M)$.  \\
\item (spherically Calabi-Yau case)  We call $(M, \omega)$ spherically Calabi-Yau, if the minimal Chern number $N$ is zero.  
\enen

Let $H: {\mathbb R}/{\mathbb Z} \times M \to {\mathbb R}$ be a smooth function.  
Set $H_t (p)=H(t,p)$.  We denote by $X_{H_t}$ the Hamiltonian vector field of $H_t$.   
We call $\ell:{\mathbb R}/{\mathbb Z} \to M$ a periodic solution of $\{X_{H_t}\}$, if $\ell$ satisfies 
$$
\frac{d}{dt} \ell(t) = X_{H_t}(\ell(t)).  
$$
We assume that  all contractible $1$-periodic solutions of $\{X_{H_t} \}$ are non-degenerate.  
Denote by ${\mathcal P}(H)$ the set of contractible 1-periodic solutions of $\{X_{H_t} \}$.  

Pick a generic $t$-dependent almost complex structure $J$ compatible with $\omega$.    
Floer chain complex $(CF_*(H, J), \delta)$ is constructed for monotone symplectic manifolds in  \cite{Floer} and for weakly monotone case in  \cite{Hofer-Salamon},  \cite{Ono}.      

Let ${\mathcal L}(M)$ be the space of contractible loops in $M$.  
Consider the set of pairs $(\ell, w)$, where $\ell: {\mathbb R}/{\mathbb Z} \to M$ is a loop and 
$w:D^2 \to M$ is a bounding disk of the loop $\ell$. 
We set an equivalence relation $\sim$ by $(\ell, w) \sim (\ell', w')$ if and only if 
$\ell = \ell'$ and  
$$
\langle [\omega], w \#(-w') \rangle = 0 \ \  \text{\rm \ and \ } \  \langle c_1(M), w \# (-w') \rangle = 0, 
$$
where $w \# (-w')$ is a spherical 2-cycle obtained by gluing 
$w$ and $w'$ with orientation reversed along the boundaries.

Then the space $\overline{\mathcal L}(M)$ of equivalence classes 
$[\ell,w]$ is a covering space of ${\mathcal L}(M)$.   
Denote by $\Pi: \overline{\mathcal L}(M) \to {\mathcal L}(M)$ 
the covering projection, and by $\Gamma$ the group of the deck
transformations of this covering, so that we have
$$
\Gamma 
=
\pi_2(M)/(\ker [\omega] \cap \ker c_1(M)).
$$
We have the weight homomorphism 
$$\int \omega : \pi_2(M) \to {\mathbb R},$$
and the corresponding Novikov ring
$\Lambda^{\mathbb Z}_{(M,\omega)} = \zz((\G))$. 


We define the action functional ${\mathcal A}_H : \overline{\mathcal L}(M) \to {\mathbb R}$ by 
$$
{\mathcal A}_H([\ell, w])= \int_{D^2} w^* \omega + \int_0^1 H(t, \ell(t)) dt.
$$  
Then the critical point set ${\rm Crit} {\mathcal A}_H$ is equal to ${\Pi^{-1}({\mathcal P}(H)})$.  
For each pair $(\ell, w)$ of $\ell \in {\mathcal P}(H)$ and its bounding disk $w$, we have the Conley-Zehnder index 
$\mu_{CZ}(\ell, w) \in {\mathbb Z}$.   
$$\mu_{CZ}: {\rm Crit} {\mathcal A}_H \to {\mathbb Z}.$$

We define $CF_k (H,J)$ by the downward completion of the free module generated by $[\ell, w] \in {\rm Crit} {\mathcal A}_H$ 
with $\mu_{CZ}([\ell,w]) = k \in {\mathbb Z}$ 
in the spirit of Novikov complex using the filtration by ${\mathcal A}_H$.  
Pick and fix a lift $[\ell, w_{\ell}]$ for each $\ell \in {\mathcal P}(H)$.  Then $CF_*(H,J)$ is a free module generated by 
$[\ell, w_{\ell}]$ over the Novikov ring 
$\Lambda^{\mathbb Z}_{(M,\omega)}$.

The boundary operator $\partial:CF_k(H,J) \to CF_{k-1}(H,J)$ is defined by counting Floer connecting orbits.  

Let $[\ell^{\pm},w^{\pm}] \in {\rm Crit} {\mathcal A}_H$.  
We denote by $\widetilde{\mathcal M}([\ell^-,w^-],[\ell^+,w^+])$ the space of the solutions $u :{\mathbb R} \times {\mathbb R}/{\mathbb Z} \to M$ 
satisfying 

\begin{equation}\label{connorb}
\frac{\partial u}{\partial \tau} + J(u(\tau, t)) ( \frac{\partial u}{\partial t} - X_{H_t}(u(\tau, t))) = 0 
\end{equation}
\begin{equation}\label{taulimit}
\lim_{\tau \to \pm \infty} u(\tau, t)= \ell^{\pm} (t) 
\end{equation}
and 
\begin{equation}\label{addcond}
[\ell^+, w^+] = [\ell^+, w^- \# u].   
\end{equation}

The group ${\mathbb R}$ acts on $\widetilde{\mathcal M}([\ell^-,w^-],[\ell^+,w^+])$ by shifting the parametrization in  $\tau$-coordinate.  
We denote by ${\mathcal M}([\ell^-,w^-],[\ell^+,w^+])$ the quotient space of $\widetilde{\mathcal M}([\ell^-,w^-],[\ell^+,w^+])$ by 
the ${\mathbb R}$-action.  Note that the ${\mathbb R}$-action is free unless $[\ell^+,w^+]=[\ell^-,w^-]$.    
We have 
$$\dim {\mathcal M}([\ell^-,w^-],[\ell^+,w^+]) = \mu_{CZ}([\ell^+,w^+]) -\mu_{CZ}([\ell^-,w^-]) -1.$$
The moduli spaces ${\mathcal M}([\ell^-,w^-],[\ell^+,w^+])$ are oriented in a compatible way, see  \cite{Floer} section 2e, \cite{FO99} section 21 and 
\cite{Flux} section 5.  
For $[\ell^{\pm},w^{\pm}]$ such that $\mu_{CZ}([\ell^+,w^+]) - \mu_{CZ}([\ell^-,w^-]) = 1$, 
${\mathcal M}([\ell^-,w^-],[\ell^+,w^+])$ is a $0$-dimensional compact oriented manifold.  
We denote by $n([\ell^-,w^-],[\ell^+,w^+])$ the order of ${\mathcal M}([\ell^-,w^-],[\ell^+,w^+])$ counted with signs.  

For $[\ell^+, w^+] \in {\rm Crit} {\mathcal A}_H$, we define 
$$
\partial [\ell^+, w^+] = \sum n([\ell^-,w^-],[\ell^+,w^+]) [\ell^-,w^-],
$$
where the summation is taken over $[\ell^-,w^-]$ such that $\mu_{CZ}([\ell^-, w^-])=\mu_{CZ}([\ell^+,w^+])-1$.  

In \cite{Hofer-Salamon}, \cite{Ono}, the Floer complex is constructed over the Novikov ring 
$\Lambda^{{\mathbb Z}/2{\mathbb Z}}_{(M,\omega)} \cong \Lambda_{(M,\omega)}^{\mathbb Z} \otimes {\mathbb Z}/2{\mathbb Z}$ with 
${\mathbb Z}/2{\mathbb Z}$-coefficients.   In order to construct the Floer complex over the Novikov ring with ${\mathbb Z}$-coefficients, 
we need an appropriate coherent system of orientations on the moduli spaces ${\mathcal M}([\ell^-,w^-],[\ell^+,w^+])$.  
Taking \cite{Floer} section 2e, \cite{FO99} section 21 into account, the argument in \cite{Hofer-Salamon} section 5 derives 
the well-definedness of the boundary operator $\partial$ and the fact that $\partial \circ \partial = 0$.  

Namely, the moduli space ${\mathcal M}([\ell_1,w_1], [\ell_2,w_2])$ with $\mu_{CZ}([\ell_2,w_2])-\mu_{CZ}([\ell_1,w_1])=2$ 
is a compact oriented $1$-dimensional manifold such that boundary is the union of the direct product 
${\mathcal M}([\ell_1w_1],[\ell,w]) \times {\mathcal M}([\ell,w],[\ell_2,w_2])$ over $[\ell, w]$ with 
$\mu_{CZ}([\ell,w])=\mu_{CZ}([\ell^+,w^+])=1$.  
This implies that the summation of 
$n([\ell_1, w_1],[\ell,w]) \cdot n([\ell,w],[\ell_2,w_2])=0$, 
hence the coefficient of 
$[\ell_1,w_1]$ in $\partial \circ \partial ([\ell_2,w_2])$ vanishes.

Hence we have 

\beth\lb{floercomplex}
Let $(M,\omega)$ be a closed weakly monotone symplectic manifold.  
For a non-degenerate $1$-periodic Hamiltonian function $H$ and a generic almost complex structure compatible with $\omega$, 
$(CF_*(H,J), \partial)$ is a ${\mathbb Z}$-graded chain complex over $\Lambda^{\mathbb Z}_{(M,\omega)}$ with integer coefficients.   
\enth 
We denote by $HF_*(H,J)$ the homology of $(CF_*(H,J), \partial)$.  

Let $H_{\alpha}, H_{\beta}$ be non-degenerate $1$-periodic Hamiltonians and $J_{\alpha}, J_{\beta}$ generic almost 
complex structures compatible with $\omega$.  
Pick a one-parameter family of smooth functions ${\mathcal H}=\{ H^{\tau} \}$ on ${\mathbb R}/{\mathbb Z} \times M$ and a one-parameter family 
${\mathcal J}=\{J^{\tau} \}$ of almost complex structures compatible with $\omega$ such that 
$H^{\tau} = H_{\alpha}$ and $J^{\tau} =J_{\alpha}$ for sufficiently negative $\tau$ and 
$H^{\tau}= H_{\beta}$ and $J^{\tau}=J_{\beta}$ for sufficiently positive $\tau$.

\beth\lb{chainequiv} Let $H_{\alpha}, H_{\beta}$ and $J_{\alpha}, J_{\beta}$
be as above.   
Then there exists a chain homotopy equivalence 
$$\Phi_{{\mathcal H},{\mathcal J}} : CF_*(H_{\alpha}, J_{\alpha}) 
\to CF_*(H_{\beta},J_{\beta}).$$
\enth 

The chain homomorphism $\Phi_{{\mathcal H}, {\mathcal J}}$ is constructed by counting {\it isolated} solutions 
$u:{\mathbb R} \times {\mathbb R}/{\mathbb Z} \to M$  joining $[\ell^-,w^-] \in {\rm Crit} {\mathcal A}_{H_{\alpha}}$ and $[\ell^+,w^+] 
\in {\rm Crit} {\mathcal A}_{H_{\beta}}$ of the following equation.  

\begin{equation}\label{chainhomeq}
\frac{\partial u}{\partial \tau} + J^{\tau}(u(\tau, t)) 
\Big( \frac{\partial u}{\partial t} - X_{H^{\tau}_t}(u(\tau, t))\Big) = 0.  
\end{equation}

For two choices $({\mathcal H}_1, {\mathcal J}_1)$ and $({\mathcal H}_2, {\mathcal J}_2)$, the chain homomorphisms 
$\Phi_{{\mathcal H}_1,{\mathcal J}_1}$ and $\Phi_{{\mathcal H}_2,{\mathcal J}_2}$ are chain homotopic.   
To construct a chain homotopy between them, we pick  homotopies $\{{\mathcal H}_s\}$, resp. $\{{\mathcal J}_s\}$, $s \in [0,1]$ 
between ${\mathcal H}_1$ and ${\mathcal H}_2$, resp. ${\mathcal J}_1$ and ${\mathcal J}_2$ and count 
{\it isolated} solutions of Equation \eqref{chainhomeq} with ${\mathcal J}_s, {\mathcal H}_s$ for some $s \in [0,1]$.    

\beth\lb{comparisonwithMorse} (\cite{Piunikhin-Salamon-Schwarz})  
Let $(M,\omega)$ be a closed symplectic manifold and $f$  a  Morse function $f$ on $M$.  
For a non-degenerate $1$-periodic Hamiltonian $H$ and a generic almost complex structure $J$ compatible with $\omega$, 
the Floer complex $(CF_*(H,J),\partial)$ is chain equivalent to the Morse complex 
$(CM_{*+n}(f) \otimes \Lambda^{\mathbb Z}_{(M,\omega)}, \partial^{\rm Morse})$.  
\enth 
For the comparison of orientation of the moduli space of solutions of Equation \eqref{connorb} and the moduli space of 
Morse gradient flow lines, see \cite{FO99} section 21.  

\bere
If $(M,\omega)$ is either monotone, spherically Calabi-Yau or the minimal Chern number $N \geq n$, 
we have a chain homotopy equivalence between $(CF_*(H,J),\partial)$ and 
$(CF_*(f,J), \partial)$ for a sufficiently small Morse function $f$.  The latter is isomorphic to 
the Morse complex $(CM_{*+n}(f) \otimes \Lambda^{\mathbb Z}_{(M,\omega)}, \partial^{\rm Morse})$, 
see \cite{Hofer-Salamon} Proposition 7.4.  

In \cite{Ono}, we introduced modified Floer homology $\widehat{HF}_*(H,J)$, which is computed in the case that $(M,\omega)$ is 
a closed weakly monotone symplectic manifold and show that 
$$\widehat{HF}_*(H,J)  \cong H_{*+n}(M;\Lambda^{\mathbb Z}_{(M,\omega)}).$$ 
(In order to work over integer coefficients, we use the orientation of the moduli space of solutions of 
Equation \eqref{connorb}, \eqref{chainhomeq} as in \cite{FO99} section 21.)   
In the end of section 6.3 \cite{FOOO}, we have an isomorphism between $\widehat{HF}_*(H,J)$ and $HF_*(H,J)$, see also 
\cite{Cieliebak-Frauenfelder} Remark 4, which also yields 
$$HF_*(H,J)  \cong H_{*+n}(M;\Lambda^{\mathbb Z}_{(M,\omega)}).$$ 
\enre

If the minimal Chern number of $(M,\omega)$ is $N$, the Floer chain complex $(CF_*(H,J),\delta)$ is $2N$-periodic, i.e., 
$$(CF_{*}, \partial) \cong (CF_{*+2N}, \partial).$$ 
(Pick an element $A \in \pi_2 (M)$ such that $\langle c_1(M), A \rangle = N$.  
Then the action of $[A] \in \pi_2(M)/(\ker [\omega] \cap \ker c_1(M))$ induces such an isomorphism of chain complexes.)  
The ${\mathbb Z}/2N{\mathbb Z}$-graded version of Floer complex $(CF_*,\delta)$ 
is a free finitely generated chain complex 
over the 
smaller Novikov ring. Namely, put
$$\Gamma_0=
\ker c_1(M)/(\ker [\omega] \cap \ker c_1(M)).$$
endow it with the \ho~ $\int \omega:\G_0\to\rr$
and consider the corresponding Novikov completion
$\Lambda^{(0){\mathbb Z}}_{(M,\omega)}=\zz((\G_0))$. 
We may also denote this ring by $\wh\L$
(observe that $\int \omega:\G_0\to\rr$ is a monomorphism).
We have the following 

\beth    
Let $(M,\omega)$ be a closed weakly monotone symplectic 
manifold with minimal Chern number $N$ and 
$H$ a non-degenerate $1$-periodic Hamiltonian on $M$.  

\noindent
(1) \ \ Then there exists a ${\mathbb Z}/2N{\mathbb Z}$-graded chain complex $(CF^{\circ}_*(H,J),\partial)$, which is freely generated by 
$\{[\ell, w_{\ell}] \vert \ell \in {\mathcal P}(H) \}$ over
$\Lambda^{(0){\mathbb Z}}_{(M,\omega)}$.

\noindent
(2) \ \ 
$(CF^{\circ}_*, \partial)$ is chain equivalent 
to the Morse complex $(CM^{\circ}_{*+n}, \partial^{\rm Morse})$ 
with the grading modulo $2N$ 
as chain complexes over  
$\Lambda^{(0){\mathbb Z}}_{(M,\omega)}$.
\enth

We can also construct Floer complex with coefficients in a local system on $M$, see \cite{Flux} section 6, and   
prove corresponding results in the  case with coefficients in a local system on $M$.    

\beth
Let $(M,\omega)$ be a closed weakly monotone symplectic manifold 
and $\rho: \pi_1(M) \to GL(r, \ff)$.  
For a non-degenerate $1$-periodic Hamiltonian $H$ and a generic almost complex structure compatible with $\omega$, 
$$HF^{\circ}_*(H,J;\rho) \cong H^{\circ}_{*+n}(M;\rho) \otimes \Lambda_{(M,\omega)}^{\ff}.$$ 
\enth

\subsection{Floer complex over a regular cover}
\label{su:on-covering}

Let $pr:\widetilde{M}\to M$ be a regular  covering  of $M$
with the covering transformation group $G$.  
Let $H$ be a non-degenerate $1$-periodic Hamiltonian $M$ and $J$ a generic $t$-dependent almost complex structures 
compatible with $\omega$.  
Denote by $\widetilde{H}$, resp. $\widetilde{J}$  the pull-back of $H$, 
resp. $J$, to ${\mathbb R}/{\mathbb Z} \times \widetilde{M}$.  
The Floer complex $(CF_*(\widetilde{H}, \widetilde{J}), \partial)$ is constructed in the spirit of  \cite{LO}.  


We will now define $\overline{\mathcal L}(\widetilde{M})$ 
similarly to $\overline{\mathcal L}(M)$. 
Consider the set of all 
pairs $(\gamma, w)$ 
where $\gamma$ is a loop 
 in $\widetilde{M}$ and 
 $w$ is  a bounding disk for $pr \circ \gamma$. Introduce in this set  
the following equivalence relation:
$(\gamma, w)$ and $(\gamma',w')$ are equivalent if 
the values of both  cohomology classes   $[\omega]$ and $c_1(M)$
on the singular sphere $w \# (-w')$ are the same. 
The set $\overline{\mathcal L}(\widetilde{M})$ of the equivalence classes is a covering space of 
$\mathcal L(\wi M)$, and the deck transformation group of the covering
is isomorphic to 
$$
\Gamma 
=
\pi_2(M)/(\ker [\omega] \cap \ker c_1(M)).
$$
The action functional  
$${\mathcal A}_{\widetilde{H}}:\overline{\mathcal L}(\widetilde{M}) \to\RRR$$
 is defined by the same formula as before, namely 
$${\mathcal A}_{\widetilde{H}}(\gamma, w)={\mathcal A}_{H}(pr \circ \gamma, w).$$

We define $CF_*(\widetilde{H}, \widetilde{J})$ by the downward completion of free abelian group generated by ${\rm Crit} {\mathcal A}_{\widetilde{H}}$ 
with respect to the action functional ${\mathcal A}_{\widetilde{H}}$.    
Pick and fix a lift $\widetilde{\ell}$ of $\ell \in {\mathcal P}(H)$ to a $1$-periodic solution of $X_{\widetilde{H}_t}$ on $\widetilde{M}$.  
Note that $pr^{-1} ({\mathcal P}(H)) = \{ g \cdot \widetilde{\ell} \vert \ell \in {\mathcal P}(H), \ g \in G \}$.  
We have $${\rm Crit}{\mathcal A}_{\widetilde{H}}= \{[g\cdot \widetilde{\ell},w] \vert [\ell, w] \in {\rm Crit}{\mathcal A}_H, \ g \in G \}.$$ 
We also pick and fix a bounding disk $w_{\ell}$ for each $\ell \in {\mathcal P}(H)$.  
Then we find that $CF_*(\widetilde{H}, \widetilde{J})$ is isomorphic to a free module generated by 
$\{ [\widetilde{\ell}, w_{\ell}]  \vert \ell \in {\mathcal P}(H) \}$ 
over 
$\Lambda^{\mathbb Z[G]}_{(M,\omega)}$.

The boundary operator is defined by counting certain isolated solutions of Equation \eqref{connorb} as follows.  
Let $[\gamma^{\pm}, w^{\pm}] \in {\rm Crit}_{\widetilde{H}}$.    
We consider the moduli space ${\mathcal M}([\gamma^-, w^-],[\gamma^+,w^+])$ of solutions of Equation \eqref{connorb} satisfying 
Condition \eqref{taulimit} with $\ell^{\pm} = pr \circ \gamma^{\pm}$, Condition \eqref{addcond} and 
that $u(\tau, 0) :{\mathbb R} \to \widetilde{M}$ lifts to a path joining $\gamma^-(0)$ and $\gamma^+(0)$. 
We set $n([\gamma^-,w^-],[\gamma^+,w^+])$ the signed count of isolated solutions in  ${\mathcal M}([\gamma^-, w^-],[\gamma^+,w^+])$.  
Since isolated connecting orbits in ${\mathcal M}([pr \circ \gamma^-, w^-],[pr \circ \gamma^+,w^+])$ are at most 
finitely many ${\mathcal M}([\gamma^-, w^-],[g \cdot \gamma^+,w^+])$ contains isolated connecting orbits.  
In other words, for fixed $[\gamma^-, w^-]$, $[\gamma^+, w^+]$, there are at most finitely many $g \in G$ such that 
$n([\gamma^-,w^-],[g \cdot \gamma^+,w^+]) \neq 0$.   
The boundary operator $\partial$ on $CF_*(\widetilde{H}, \widetilde{J})$ is given by 
$$
\partial [\gamma^+, w^+] = \sum n([\gamma^-,w^-],[\gamma^+,w^+]) [\gamma^-,w^-].
$$
It is clear that $\partial$ is linear over 
$\Lambda^{\mathbb Z[G]}_{(M,\omega)}$.

Keeping attention on the homotopy classes of paths $u(\tau, 0):{\mathbb R} \to M$ of solutions 
of Equation \eqref{connorb}, \eqref{chainhomeq}, 
the proofs of Theorems \ref{floercomplex}, \ref{chainequiv} and \ref{comparisonwithMorse} works for 
the case of $CF_*(\widetilde{H},\widetilde{J})$.  
For example, we show the fact that  $\partial \circ \partial =0$ in the following way.  
In subsection \ref{su:review},  we recalled that $(CF_*(H,J), \partial)$ is a chain complex using the moduli space 
${\mathcal M}([\ell_1,w_1],[\ell_2,w_2])$ with $\mu_{CZ}([\ell_2,w_2]) -\mu_{CZ}([\ell_1,w_1])=2$.  
The projection $pr:\widetilde{M} \to M$ gives an identification of 
${\mathcal M}([\gamma_1,w_1],[\gamma_2,w_2])$ and the subspace of 
${\mathcal M}(pr \circ \gamma^-,w^-,pr \circ \gamma^+,w^+])$ consisting of connecting orbits $u$ such that 
$u(\tau,0)$ lifts to a path joining $\gamma_1(0)$ and $\gamma_2(0)$.  
The boundary of this subspace is the union of the direct product 
${\mathcal M}([\gamma_1,w_1],[\gamma,w]) \times {\mathcal M}([\gamma,w],[\gamma_2,w_2])$ 
such that $\mu_{CZ}([\gamma,w])=\mu_{CZ}([\gamma^+,w^+])-1$, which is 
identified with the union of the space of pairs $(u_1, u_2)$ of 
${\mathcal M}([pr \circ \gamma_1,w_1],[pr \circ \gamma,w]) \times {\mathcal M}([\pr \circ \gamma,w],[pr \circ \gamma_2,w_2])$ 
such that the concatenation of the paths 
$u_1(\tau, 0)$ and $u_2(\tau, 0)$ lifts to a path joining 
$\gamma_1(0)$ and $\gamma_2(0)$.  
Hence, by looking at the components of 
${\mathcal M}(pr \circ \gamma_1,w_1,pr \circ \gamma_2,w_2])$ 
with $u(\tau, 0)$ in the prescribed homotopy class of paths joining $pr \circ \gamma_1(0)$ and $pr \circ \gamma_2(0)$, 
we find that $CF_*(\widetilde{H},\widetilde{J})$ is a chain complex.  
This Floer complex is periodic with respect to the degree shift by $2N$.  Hence we can also obtain 
${\mathbb Z}/2N{\mathbb Z}$-graded chain complex, which we denote by 
$CF^{\circ}_*(H,J)$, $CF^{\circ}_*(\widetilde{H}, \widetilde{J})$, etc.

\beth\label{t:w-monoton-cov}
Let $\widetilde{M}\to M$ be a regular  covering of a closed 
weakly monotone symplectic manifold $(M, \omega)$ 
with minimal Chern number $N$. 
Let $G$ be the structure group 
of the covering.

\noindent
(1) \ \ For a non-degenerate $1$-periodic Hamiltonian $H$ on $M$, there exists a ${\mathbb Z}/2N{\mathbb Z}$-graded chain complex 
$(CF^{\circ}_*(\widetilde{H}, \widetilde{J}), \partial)$ such that $CF^{\circ}_*(\widetilde{H},\widetilde{J})$ is a free module generated by 
$\{ [\widetilde{\ell}, w_{\ell}] \vert \ell \in {\mathcal P}(H) \}$ over 
$\Lambda^{(0) {\mathbb Z}[G]}_{(M,\omega)}$.

\noindent
(2) \ \ 
Let $f$ be a Morse function on $M$.  
Then for a non-degenerate $1$-periodic Hamiltonian $H$ and a generic almost complex structure $J$ compatible with $\omega$, 
$(CF^{\circ}_*(\widetilde{H},\widetilde{J}), \partial)$ is chain equivalent to 
the Morse complex 
$(CM^{\circ}_{*+n}(f \circ pr) \otimes_{\mathbb Z[G]} \Lambda^{(0){\mathbb Z}[G]}_{(M,\omega)},
\partial^{\rm Morse})$ 
of the Morse function $f\circ pr: \wi M\to\rr$
 with coefficients in 
 $\Lambda^{(0){\mathbb Z}[G]}_{(M,\omega)}$.
\enth

In the case of arbitrary closed symplectic manifolds
we have an analog of this result over the field $\qq$.

\beth\label{t:gen-case-cover}
Let $\widetilde{M}\to M$ be a regular  covering of a closed 
symplectic manifold $(M, \omega)$ 
with minimal Chern number $N$. 
Let $G$ be the structure group 
of the covering.

\noindent
(1) \ \ For a non-degenerate $1$-periodic Hamiltonian $H$ on $M$, there exists a ${\mathbb Z}/2N{\mathbb Z}$-graded chain complex 
$(CF^{\circ}_*(\widetilde{H}, \widetilde{J}), \partial)$ such that $CF^{\circ}_*(\widetilde{H},\widetilde{J})$ is a free module generated by 
$\{ [\widetilde{\ell}, w_{\ell}] \vert \ell \in {\mathcal P}(H) \}$ over 
$\Lambda^{(0){\mathbb Q}[G]}_{(M,\omega)}$.   

\noindent
(2) \ \ 
Let $f$ be a Morse function on $M$.  
Then for a non-degenerate $1$-periodic Hamiltonian $H$ and a generic almost complex structure $J$ compatible with $\omega$, 
$(CF^{\circ}_*(\widetilde{H},\widetilde{J}), \partial)$ is chain equivalent to 
the Morse complex 
$(CM^{\circ}_{*+n}(f \circ pr) \otimes_{\mathbb Q[G]} 
\Lambda^{(0) \mathbb Q[G]}_{(M,\omega)}, \partial^{\rm Morse})$ 
of the Morse function $f\circ pr: \wi M\to\rr$
 with coefficients in 
$\Lambda^{(0){\mathbb Q}[G]}_{(M,\omega)}$.  
\enth

\section{Invariants of chain complexes: $\zz$-graded case}
\label{s:z-graded}

The sections 
\ref{s:z-graded} and 
\ref{s:z2-graded}
are purely algebraic. 
We introduce some invariants of chain complexes, which will
be applied in Section \ref{s:estimates}
to obtain lower bounds for 
$p(H)$.

\subsection{Definition of invariants $\mu_i$}${}$
\label{su:def-mu}

Recall that a ring $R$ is called {\it an IBN-ring}, if the cardinality
of a base of a free $R$-module does not depend on the 
choice of the base.
Any principal ideal domain (PID) is an IBN-ring.
The group ring of any group with coefficients
in a PID is an IBN-ring. All the rings which we consider 
in this paper will be IBN-rings. 

\bede\label{d:mu}
For a free based finitely generated module $A$ over 
an IBN-ring $R$ we denote
by $m(A)$ the cardinality of any  base of $A$.
It will be called {\it the rank of $A$}.

Let $C_*=\{C_n\}_{n\in\zz}$ be a \ffg~ \cco~ over a ring $R$.
Denote by $m_i(C_*)$ the number $m(C_i)$.
The minimum of the numbers $m_i(D_*)$, where $D_*$ ranges over
the set of all \fbfg~\cco es chain equivalent to $C_*$,
will be denoted by $\mu_i(C_*)$.
\end{defi}
Observe that the chain complexes which we consider
are not supposed to vanish  in 
negative degrees.

Our aim in this section is to develop efficient tools for computing
the invariants $\mu_i(C_*)$ for the case of chain complexes
arising in the applications to the Arnold conjecture.

We will use  here  the terminology from 
Subsection \ref{su:def-nov}. Namely $G$ is a group, 
$T$ is a free abelian finitely generated group, 
$\xi:T\to\rr$ is a monomorphism. We consider the rings 
$$\LL=\zz[G]((T)),
\ \ 
\wh\L = \zz((T)), 
\ \ 
\FF = \ff((T)).
$$
The ring $\LL$ is an IBN-ring since it has an epimorphism onto
$\wh\L$. Similarly, $\LL_\ff$ is an IBN-ring.

\noindent 
We will use the following notation
throughout the rest of the paper:
\begin{quote}
$X$ is a  connected finite CW-complex,
$\wi X\to X$ is a regular covering with a structure group
$G$, so that we have an epimorphism 
$\pi_1(X)\to G$, 

\begin{equation}\lb{f:c-x}
\CCCC_*(X)=C_*(\wi X)\tens{\zz G}\zz[G]((T)).
\end{equation}
\end{quote}
\noindent
Then $\CCCC_*(X)$ is a free finitely 
generated chain complex over $\LL=\zz[G]((T))$. 
Observe the following isomorphism of $\LL$-modules:
\begin{equation}\lb{f:h-zero}
 H_0(\CCCC_*(X)) \approx \wh\L.
 \end{equation}

\subsection{Lower bounds provided by the 
cohomology with local coefficients}
\label{su:loc-coef}
${}$

Let $\r:G\to\GL(r,\ff)$ be a representation. 
We denote by $b_i(X,\r)$ the Betti numbers
of $X$ \wrt~ the local coefficient system 
induced by $\r$.
Put
\begin{equation}\lb{f:betti}
 \b_i(X,\r)= \frac 1rb_i(X,\r).
\end{equation}
We will need the following basic Lemma.
\bele\lb{l:betti}
Let $D_*$ be a free finitely generated 
chain complex over $\LL$, chain equivalent to 
$\CCCC_*(X)$. 
Let $\r:G\to\GL(r,\ff)$ be a representation. 
Then there is a chain complex $E_*$ over $\FF$ such that 
\been\item$\dim_\FF E_k = r\cdot m_k(D_*)$,
\item $\dim_\FF H_k(E_*)=b_k(X,\r)$.
\enen
\enle
\Prf 
The vector space $\ff^r$ has the structure of $\zz[G]$-module 
via the representation $\r$. This structure induces 
a natural structure $\wh\r$ of $\LL$-module on
$\FF^r$ by the following formula
$$
(\sum_i a_i g_i ) (\sum_j v_j h_j) 
=
\sum_{i,j} a_i(v_j) g_i h_j\ \ \ 
{\rm with ~ ~ } a_i\in \zz[G], \  \ v_j\in \ff^r, \  g_i, \  h_j \in T.
$$
\noindent
We have also the representation 
$\r_0:G\to \GL(r, \FF)$ 
obtained as the composition 
of $\r$ with the embedding $\ff\arrinto \FF$.
Put $E_*=D_*\tens{\wh\r} \FF^r$. 
Then $\dim_\FF E_k = r\cdot m_k(D_*)$, and
$$
\CCCC_*(X)\tens{\wh\r}\FF^r = 
\Big(C_*(\wi X)\tens{\zz[G]}\LL\Big)\tens{\wh\r}\FF^r
\approx
C_*(\wi X)\tens{\r_0}\FF^r
\approx 
\Big(C_*(\wi X)\tens{\r} \ff^r\Big)\tens{\ff}\FF
$$
so that
$$
\dim_\FF H_k(E_*) = \dim_\FF H_k\Big(\CCCC_*(X)\tens{\wh\r} \FF^r\Big) 
= b_k(X,\r).\qquad \qs 
$$
\noindent
The next proposition follows.
\bepr\lb{p:bn}
We have
$$
\mu_i(\CCCC_*(X))
\geq
\b_i(X,\r).
$$
\enpr
\Prf
Let $D_*$ be any chain complex over $\LL$ which is 
chain equivalent to $\CCCC_*(X)$.
Pick a chain complex 
$E_*$ constructed in the previous Lemma.
We have 
$$
r\cdot m_k(D_*) = \dim_\FF E_k\geq 
\dim_\FF H_k(E_*)=b_k(X,\r) = r\cdot \b_i(X,\r). \qquad \qs
$$


For $i=1$ we have a slightly stronger version of this estimate.
In order to prove it, we need a lemma.
\bele\lb{l:splittingoff}
Let $A_*$ be a \fbfg~ \cco~ over a ring $R$,
such that $A_*=0$ for $*\leq l-2$ and $H_{l-1}(A_*)=0$.
Then there exists a chain complex $B_*$ such that
\been\item
$B_j\approx A_j$\ \ for\ $j\geq l+2$,\  and \ $j=l$.
\item
$B_{l+1}=A_{l+1}\oplus A_{l-1}$,
\item
$B_j=0$\ for\ $j\leq l-1$.
\enen
\enle
\Prf
Let $T_*$ be the chain complex
$$\{0\leftarrow A_{l-1}\arrl {Id} A_{l-1}\leftarrow 0 \}$$
concentrated in degrees $l$ and $l+1$.
By the Thickening Lemma 
(\cite{Sharko}, Lemma 3.6, p. 56)
the chain complex $C_*=A_*\oplus T_*$
is isomorphic to the chain complex
$$
C'_*
=
\{ 0\leftarrow A_{l-1}\arrl {\pr_l} A_{l-1}\oplus A_{l}
\arrl {\pr_{l+1}} A_{l-1}\oplus A_{l+1}
\leftarrow A_{l+2}\la \ldots \}
$$
with $\pr_l(x,y)=x$. Splitting off the chain complex
$$\{ 0\leftarrow A_{l-1}\arrl {Id} A_{l-1}\leftarrow 0 \}$$
concentrated in degrees $l-1$ and $l$, 
we obtain the required chain complex $B_*$. $\qs$

\bepr\lb{p:bn1}
We have
$$
\mu_1(\CCCC_*(X))
\geq
\b_1(X,\r)+1 
$$
for any representation $\r$ such that $H_0(X,\r)=0$.
\enpr
\Prf
Let $D_*$ be any \fbfg~ \cco~ over $\LL$, 
chain equivalent to $\CCCC_*(X)$.
An easy  induction argument 
using Lemma \ref{l:splittingoff}
shows that $D_*$ is chain 
equivalent to
a \fbfg~ \cco~ $D_*'$ such that $D_i'=0$ for $i\leq -2$ and
$D'_i=D_i$ for $i \geq 1 $.
Put 
$$\a=m (D'_{-1}), 
\ \
\b=m (D'_{0}),
\ \
\g=m (D'_{1}).
$$
The homology of $D'_*$ and 
of $D_*\tens{\wh\r}\wh\L^r_\ff$
vanishes in degree $-1$.
Applying Lemma \ref{l:betti}
to the trivial 1-dimensional representation $\r_0$
we find a chain complex 
$$E_*
=
\{ 0\arl
 \FF^{\a }
\arrl {\pr_0} 
 \FF^{\b }
\arrl {\pr_1}
 \FF^{\g }
\arl \ldots \ \ \}\ \ 
$$
such that $H_{-1}(E_*)=0$ and $H_0(E_*)\approx \FF$;
this implies $\b\geq \a+1$.
Applying Lemma \ref{l:betti}
to the  representation $\r$
we obtain  a chain complex
$$E'_*
=
\{ 0\arl
 \FF^{\a r}
\arrl {\pr'_0} 
 \FF^{\b r}
\arrl {\pr'_1}
 \FF^{\g r}
\arl \ldots \ \  \} \ \ 
$$
with $\dim_\FF\Ker\pr'_1 \geq b_1(X,\r)$. Since $H_0(X,\r)$
and $H_0(E'_*)$ vanish, we obtain
$$
r\g\geq b_1(X,\r)+r(\b-\a) \geq b_1(X,\r)+r.\qquad \qs
$$ 
\bere\label{r:irred-null}
Observe that the condition $H_0(X,\r)=0$ is true for every 
non-trivial irreducible representation $\r$.
\end{rema}
When $\pi_1(X)$ is a perfect group the above methods allow to 
obtain a lower bound for $\mu_2(X)$:
\bepr\lb{p:mu-2}
Assume that the covering $\wi X\to X$ is the universal covering of $X$,
so that in particular $\pi_1(X)\approx G$. Assume that $G$ is a perfect
finite group. Then $\mu_2(X)\geq b_2(X)+2$.
\enpr
\Prf
Let $D_*$ be any \fbfg~ \cco~ over $\LL$, 
chain equivalent to $\CCCC_*(X)$.
Consider the chain complex $D'_*$ constructed in the 
proof of the Lemma 
\ref{p:bn1}. 
Applying Lemma \ref{l:betti}
to the trivial 1-dimensional representation $\r_0$
we find a chain complex 
$$E_*
=
\{ 0\arl
 \FF^{\a }
\arrl {\pr_0} 
 \FF^{\b }
\arrl {\pr_1}
 \FF^{\g }
\arrl {\pr_2}
 \FF^{\d }
\arl 
\ldots \ \ \}\ \ 
$$
such that $H_{-1}(E_*)=0$ and $H_0(E_*)\approx \FF$;
We have $\dim\Ker\pr_1= \g-(\b-\a)+1$. 
Since $G$ is perfect, we have $\dim H_1(E_*)=b_1(X)=0$, therefore
$\dim\Im\pr_2'= \g-(\b-\a)+1$, so that
$$
\d\geq b_2(X) + \g-(\b-\a)+1.
$$
Choose an irreducible representation $\r$, such that $b_1(X,\r)\geq 1$
(this is possible by Theorem \ref{t:d-mu1}). 
Applying Lemma \ref{l:betti}
to the   representation $\r$ we deduce
$\g-(\b-\a) \geq  \b_1(X,\r) >0$, so that 
$\d\geq  b_2(X) + 2.$ $\qs$

\subsection{The invariant $\mu_1$: the case of finite groups}
\label{su:finite-groups}

The results of the previous subsection
allow to obtain a complete result for the 
invariant $\mu_1$ in the case of finite groups.
For a group $G$ denote by $d(G)$ the minimal possible number of 
generators of $G$ and by $\delta(G)$ the minimal possible number
of generators of the augmentation ideal of $\zz[G]$
as a $\zz[G]$-module.
The next theorem 
is a reformulation of a well-known result in the
cohomological theory of finite groups
(see, for example, \cite{Roggenkamp}, Corollary 5.8, p. 191).

\beth\lb{t:d-mu1}
Let $G$ be a finite group. Then $\d(G)$
equals the maximum of 
two numbers $A(G)$ and $B(G)$, defined below.
$$A(G)=\underset{p,\r}{\max}\Big( \Big\ldbrack\frac 1r b_1(G,\r)\Big\rdbrack
+1\Big),$$
where the maximum is taken over all prime 
divisors $p$ of $|G|$ and all 
the irreducible non-trivial representations
$\r:G\to \GL(r,\ff_p)$.
\footnote{The symbol $\ldbrack z \rdbrack$
denotes the minimal integer $k\geq z$.}
$$B(G)=\underset{p}{\max}\ b_1(G,\ff_p),$$
where the maximum is taken 
over all prime 
divisors $p$ of $|G|$.
$\qs$
\enth

\bere\label{r:g-x}
If $\phi : G\to K$ is a group  epimorphism,
$V$ a $K$-module, and we endow $V$ with
a structure of $G$-module via $\phi$, then the induced 
\ho~ $H_1(G,V)\to H_1(K, V)$ is surjective.
\enre

Recall from \rrf{f:c-x}
that $X$ denotes a connected finite CW-complex, and 
$\wi X\to X$ is a regular covering with a structure group
$G$, so that we have an epimorphism 
$\pi_1(X)\to G$. 
We denote by $\CCCC_*(X)$
the chain complex
$C_*(\wi X)\tens{\zz G}\LL$.

\beth\lb{t:mu1}
Let $G$ be a finite group with a group epimorphism $\pi_1(X) \to G$.
Then $\mu_1(\CCCC_*(X)) \geq\delta(G)$.
\enth
\Prf
Let $D_*$ be a \fbfg~ \cco~chain equivalent to 
$\CCCC_*(X)$. The  inequality 
$\mu_1(\CCCC_*(X)) \geq \delta(G)$
follows immediately
from Propositions 
\ref{p:bn}, \ref{p:bn1} together with Theorem 
\ref{t:d-mu1}
and Remarks \ref{r:irred-null}, \ref{r:g-x}. $\qs$

The invariant $\delta(G)$ of a finite group
has the following properties:

\been\item
$\delta(G)=d(G)$ if $G$ is solvable 
(K. Gr\"unberg's theorem \cite{GruenbergS}, see also 
\cite{Roggenkamp}, Theorem 5.9).
\item 
 $\delta(G)=1$ if and only if $G$ is 
cyclic (see \cite{Roggenkamp}, Lemma 5.5).
\enen
The second point implies also that 
$\d(G)$ equals 2 for any simple non-abelian
group $G$ (since $d(G)=2$ for such a group).
\beco\lb{c:mu}
Let $G$ be a finite group with a group epimorphism $\pi_1(X) \to G$.
We have
\been\item
$\mu_1(\CCCC_*(X))\geq d(G)$ if $G$ is solvable or simple.
\item
$\mu_1(\CCCC_*(X))\geq 1$ for every non-trivial group $G$.
\item
$\mu_1(\CCCC_*(X))\geq 2$ if $G$ is not cyclic.
\enen
\enco

\subsection{The invariant $\mu_1$: the case of infinite groups}
\label{su:inf-groups}

If the group $G$ is infinite, it is more difficult
to give computable lower bounds for $\mu_1(\CCCC_*(X))$.
In this section we prove that $\mu_1(\CCCC_*(X))>0$.
Recall the ring $\LL=\zz[G]((T))$ and its subring $\wh\L=\zz((T))$.
Observe that the $\wh\L$-module $\LL$ has no torsion.
The ring $\wh\L$ is also a module over $\LL$
via the augmentation \ho~ $\ve$ (see Subsection \ref{su:def-mu}).

\bele\label{l:non-embedded}
A free $\LL$-module contains no submodule isomorphic to $\wh\L$.
\enle
\Prf
Assume that there is an embedding $i:\wh\L\to\LL^n$.
There is a projection $p:\LL^n\to \LL$ such that $p\circ i $ is non-trivial.
Since $\LL$ has no $\wh\L$-torsion, the \ho~ 
$p\circ i $ is an embedding. Put $a=(p\circ i)(1)\in\LL$.
Multiplying 
$a$ by a suitable element of $T$ if necessary, 
we can consider that 
$a=\a\cdot {\mathbf{ 1}}+\a'$
where $\a\in\zz[G]$, and $\a'$ is a power series 
in monomials $g_i\in T$ with 
$\xi(g_i)<0$.
For every $g\in G$ we have then $g\cdot a = \l\cdot a$ with 
with some $\l\in \wh\L$, which implies $g\a=l\a$ for some $l\in\zz$. The last property 
is impossible, since $G$ is infinite, and $\a$ is a finite linear
combination of elements of $G$. $\qs$

\beco\label{c:h-0}
Let $D_*$ be a \fbfg~\cco~ over $\LL$.
Assume that $H_0(D_*)\approx \wh\L$.
Then $D_1\not=0$.
\enco
\Prf
If $D_1=0$, then $H_0(D_*)$ is isomorphic to the kernel of
the boundary operator $\pr_0:D_0\to D_{-1}$, therefore
$H_0(D_*)$ is a submodule of  free $\LL$-module, which 
contradicts to Lemma \ref{l:non-embedded}. $\qs$

\noindent
The next proposition follows.
\bepr\lb{p:betti1}
If $G$ is infinite, then $\mu_1(\CCCC_*(X))\geq 1$.
\enpr

\bere\lb{r:inf-fin}
This proposition is valid also for non-trivial finite groups,
see Theorem \ref{t:mu1}.
\enre

\subsection{The invariant $\mu_1$: the general case}
\label{su:stable-num}

\bede\label{d:stablegener}
Let $R$ be a ring, and $N$ be a module over $R$.
The minimal number $s$ such that there exists
an epimorphism 
$R^{s+r}\rightarrow N\oplus R^r$
will be called {\it the stable number of generators
of $N$}, and denoted by $\s(N)$.
The stable number of generators of the augmentation ideal
$\Ker\ve:\LL\to\wh\L$ 
of the  ring $\LL$
will be denoted $\s(G)$.
\end{defi}
\bere\label{r:stable-unstable}
Using Theorem  \ref{t:d-mu1}
it is easy to show that 
$\s(G)=\d(G)$ for any finite group.
It does not seem that this equality  holds in general,
although we do not have a counter-example at present.
\enre

\bepr\label{p:stblegener}
$\mu_1(\CCCC_*(X))\geq \s(G).$
\enpr
\Prf
Let $D_*$ be any \fbfg~ \cco~ over $\LL$, 
chain equivalent to $\CCCC_*(X)$.
Similarly to Proposition \ref{p:bn1}
we construct 
a \fbfg~ \cco~ $D_*'$ chain equivalent to $D_*$
such that $D_i'=0$ for $i\leq -2$ and
$D'_i=D_i$ for $i \geq 1 $.
Denote
$ D'_{-2}$ by $A$ and $ D'_{-1}$ by $B$.
Similarly to the proof of \ref{p:bn1} 
we deduce $m(D_0)\geq m(B)-m(A)+1$.
We have $H_i(D_*)=0$ for $i<0$. Applying Lemma  
\ref{l:splittingoff}
we obtain a chain complex 
$$
D_*'=\{\ldots 0\la D'_{-2}\arrl {\pr_{-1}}
 D'_{-1}\la D_0 \la D_1 \la \ldots\ \ 
\}
$$
Applying Lemma \ref{l:splittingoff}
two more times, we obtain
a chain complex
$$
D''_*=\{\ldots 0\leftarrow D_0\oplus A\la D_1\oplus B \la D_2 \la \ldots\ \ 
\}
$$
chain equivalent to $D_*'$.
Since $H_0(D'_*)\approx \wh\L$, we have an exact  sequence 
$$
0\la\wh\L \arrl \phi D_0\oplus A \arrl \psi D_1\oplus B. 
$$
Add to it the  exact sequence
$\{ 0\la 0 \la \LL \arrl {Id} \LL \la 0 \}$.
By the Thickening Lemma 
the result is isomorphic to the following exact sequence
$$
0\la\wh\L \arrl \chi \LL\oplus D_0\oplus A 
\arrl \phi \LL\oplus  D_1\oplus B 
$$
where $\chi(f,d,a)=\varepsilon(f)$.
Let $J(G)= \Ker(\ve:\LL\to\wh\L)$. 
We have $\Ker\chi=J(G)\oplus D_0\oplus A$, 
so that $\phi$ is an epimorphism of a free $\LL$-module of rank
$m(B)+m(D_1)+1$ onto the sum of $J(G)$ 
and a free $\LL$-module of rank 
$m(A)+m(D_0)$. Thus 
$$\s(G)\leq m(B)- m(A) +1 -m(D_0) +m(D_1) \leq m(D_1). 
\qquad\qs$$

\subsection{The invariant $\mu_2$}
\label{su:mu2}

 The results about this invariant are less complete,
than for  $\mu_1$: we have two different lower bounds
for $\mu_2(\CCCC_*(X))$ (Proposition \ref{p:b1} and Corollary
\ref{c:g-finite-mu2}), none of them is optimal in general. 
Denote by $B_1(X)$ the maximum of numbers
$\b_1(X,\r)-\b_0(X,\r)$ where $\r$ ranges over all representations of 
$G$. 
\bepr\label{p:b1}
For every representation $\r:G\to \GL(r, \ff)$
we have
\begin{equation}\label{f:mu2_B1}
 \mu_2(\CCCC_*(X)) 
 \geq
 B_1(X) +\b_2(X,\r)-\b_1(X,\r)+\b_0(X,\r). 
\end{equation}
\enpr
\Prf
Let $D_*$ be any \ffg~ \cco~ over $\LL$,
chain equivalent to $\CCCC_*(X)$. Similarly to Proposition
\ref{p:bn1} we can assume that $D_i=0$ for $i\leq -2$. 
Put 
$$\a=m (D_{-1}), 
\ 
\b=m (D_{0}),
\
\g=m (D_{1}), \ 
\d = m (D_{2}). 
$$
Apply Lemma \ref{l:betti} and let $E_*$ be the corresponding chain complex. 
Denote by $Z_0$
the space of cycles of degree $0$ of this complex;
then $\dim_\FF Z_0=r(\b-\alpha)$.
Consider the chain complex 
$$0\la Z_0\la E_1 \la E_2 \la \ldots $$
of vector spaces over $\FF$. 
Its Betti numbers are 
equal to the Betti numbers of $X$ with coefficients in $\r$, 
and applying the strong Morse inequalities we obtain 
 the following:
\begin{gather}
 \b-\a\geq \b_0(X,\r); \label{f:1} \\
 \g-(\b-\a) \geq \b_1(X,\r) -\b_0(X,\r); \label{f:2}\\
 \d -\g+\b-\a \geq \b_2(X,\r) - 
 \b_1(X,\r) +\b_0(X,\r).\label{f:3}
\end{gather}
The inequality \rrf{f:2} implies that $\g-\b+\a\geq B_1(X)$. 
Now the proposition follows from \rrf{f:3}. $\qs$

\beco\label{c:mu2_first}
Assume that $G$ is finite and the epimorphism 
$\pi_1(X)\to G$ is an isomorphism. Then
$$\mu_2(\CCCC_*(X))\geq \d(G)  - b_1(X,\ff) + b_2(X,\ff).$$

\enco
\Prf It follows from Theorem \ref{t:d-mu1} 
that $B_1(X)+b_0(X,\ff)\geq \d(G)$. $\qs$

\bere\lb{r:perfect-again}
If $G$ is a finite perfect
group then $b_1(X,\ff)=0$, and $\d(G)\geq 2$;
thus we recover the Proposition \ref{p:mu-2}.
\enre

Now we will give a lower bound for $\mu_2(X)$ in terms
of a numerical  invariant depending only on $G$
and related to the invariant $D(G)$ 
(see Introduction).
Up to the end of this Section we assume that $G$ is finite and the epimorphism 
$\pi_1(X)\to G$ is an isomorphism.
In this case the natural inclusion 
$\wh\L[G]\arrinto \zz[G]((T))$ is an isomorphism.
We will make no difference between these two rings;
observe also that
$$
\CCCC_*(X) = C_*(\wi X)\tens{\zz} \wh\L.
$$

\bede\label{d:muG}
Let $R$ be a commutative ring and
$$
\FF_*=\{0\la R\la F_0\la F_1\la \ldots \}
$$
be a free $R[G]$-resolution of the trivial
$R[G]$-module $R$; put $m_i(\FF_*)=m(F_i)$.
The minimum of $m_i(\FF_*)$ over all free resolutions
of $R$ will be denoted by $\mu_i(G,R)$.
\footnote{ Our terminology here differs from that of 
the Swan's paper 
\cite{Swan64}. }
If $R=\zz$ we abbreviate $\mu_i(G,R)$  to $\mu_i(G)$.
\end{defi}
The following properties are easy to prove:
\been
\item For any ring $R$ we have 
$\mu_i(G,R)\leq \mu_i(G)$.
\item
$\mu_1(G)=\d(G)$. 
\item $D(G)\geq 
\mu_1(G)+\mu_2(G)$. 
\enen
We will now introduce a similar notion 
appearing in the context of $\zz$-graded complexes.
\bede\lb{d:z-resolution}
Let $R$ be a commutative ring.
A  $\zz$-graded \cco~  of \ffg~ $R[G]$-modules
$$
\EE_*=
\{\ldots \la \EE_{-n}\la \ldots \la \EE_0\la\ldots \EE_{n}\la \ldots \}
$$
is called 
{\it a $\zz$-graded resolution of the trivial $R[G]$-module $R$}
if 
\been\item
$H_*(\EE_*)=0$ for every $i\not=0$ and 
$H_0(\EE_*)\approx R$.
\item$\EE_{-n}=0$ for every $n\geq 0 $ sufficiently large.
\enen
The minimum of $m_i(\EE_*)$ over all  $\zz$-graded  resolutions
of $R$ will be denoted by $\bar\mu_i(G,R)$.
If $R=\zz$ we abbreviate $\bar\mu_i(G,R)$  to $\bar\mu_i(G)$.
\end{defi}
We have obviously $\bar\mu_i(G,R)\leq \mu_i(G,R)$.
The next proposition follows from the fundamental result of R. Swan 
\cite{Swan64}. 
\bepr\label{p:z-n}
We have  $\mu_i(G)=\bar\mu_i(G)$ for $i=1,2$.
\enpr
\Prf
Let $\EE_*$ be a
$\zz$-graded resolution.
Similarly to Proposition \ref{p:bn1}
we can assume that $\EE_i=0$ for $i\leq -2$;
put
$$f_0= m(\EE_0)-m(\EE_{-1}),
\ \ 
f_i= m(\EE_i) 
\ \ {\rm for } \ \ 
i \geq 1. 
$$
Let $\r:G\to \GL(r,\ff)$ be any irreducible representation.
Put $\EE_*^\r= \EE_*\tens{\r}\ff^r$,
and let $Z_0^\r$ be the vector space of cycles of degree $0$.
Then $\dim Z_0^\r = rf_0$.
By the strong Morse inequalities applied to the chain complex 
$$0\la Z_0^\r \la  \EE^\r_1 \la \EE^\r_2 \la \ldots $$
we have 
\begin{gather}
 f_0\geq \b_0(\EE_*,\r); \label{ee:1} \\
 f_1-f_0\geq \b_1(\EE_*,\r) -\b_0(\EE_*,\r); \label{ee:2}\\
 f_2-f_1+f_0 \geq \b_2(\EE_*,\r) -\b_1(\EE_*,\r) +\b_0(\EE_*,\r).\label{ee:3}
\end{gather}
By the Swan's theory
(\cite{Swan64}, Th. 5.1, Corollary 6.1, and Lemma 5.2)
there exists a free resolution $\FF_*$
of $\zz$ over $\zz[G]$ such that
$m_i(\FF_*)=f_i$ for $i=0,1,2$. 
Therefore $\mu_i(G)\leq f_i$ for $i=1,2$. 
The proposition follows.
$\qs$

\bere\label{r:swan-all-i}
The proposition is valid for all $i\geq 1 $, with some mild
restrictions on $G$ (see Theorem 5.1 of \cite{Swan64}.)
\enre

A similar method proves the next proposition.
\bepr\label{p:z-lambda}
$\bar\mu_i(G, \wh\L)=\mu_i(G)$. $\qs$
\enpr
Now we can obtain the estimate for $\mu_2(\CCCC_*(X))$.
\bepr\label{p:mu2x}
We have $\mu_2(\CCCC_*(X))\geq \bar\mu_2(G, \wh\L)$.
\enpr
\Prf
Let $D_*$be a \ffg~ \cco~ 
chain equivalent to $\CCCC_*(X)$.
Then 
$$
H_i(D_*)=0
\ \ {\rm for }  \ 
i<0,
\ \ {\rm and }  \ 
H_0(D_*)\approx \wh\L,
\ \ {\rm and }  \ 
H_1(D_*)=0.
$$
Using the standard procedure of killing the homology groups 
of a chain complex, we embed $D_*$ into a free  
chain complex $D'_*=D_*\oplus E_*$ such that
$D_*'$ is finitely generated in each dimension 
and 
$$
E_i=0
\ \ {\rm for }  \ 
i\leq 2,
\ \ {\rm and }  \
H_0(D_*')=\wh\L
\ \ {\rm and }  \
H_i(D_*')=0
\ \ {\rm for }  \
i\not=0.
$$
Then 
$m_2(D_*)=m_2(D'_*)\geq \mu_2(G). $
The proposition follows. $\qs$

The next Corollary is immediate. 
\beco\label{c:g-finite-mu2}
$\mu_2(\CCCC_*(X))\geq \mu_2(G)$. $\qs$
\enco

\section{Invariants of chain complexes: $\zz/k\zz$-graded case}
\label{s:z2-graded}

\bede\lb{d:2graded}
Let $R$ be a ring, and $k\in\nn, k\geq 2$.
A {\it $\zz/k\zz$-graded chain complex } 
is a family of 
\fbfg~ $R$-modules $A_i$ 
indexed by  $i\in\zz/k\zz$
together with \hos~ 
$\pr_i: A_{i}\to A_{i-1}$,
satisfying $\pr_i\circ\pr_{i+1} = 0.$
\end{defi}
Given $k\in\nn$ and a \fbfg~ $\zz$-graded chain complex $C_*$,
one constructs 
a $\zz/k\zz$-graded chain complex $C ^\circ_\stakr$ as follows:
$$
C^\circ_{i}=\bigoplus_{s\equiv i(k)} C_{s}.
$$
In this section we will be working with the $\zz/k\zz$-graded 
chain complex induced by $\CCCC_*(X)$
(see the definition
\rrf{f:c-x}). It will be denoted by 
$\CCCC^\circ_\stakr(X)$, where 
\begin{equation}\lb{f:svertka}
\CCCC^\circ_i(X)
=
\bigoplus_{s\equiv i(k)}\CCCC_{s}(X).
\end{equation}
\bede\lb{d:mumu}
Let $C_\stakr$ be a 
$\zz/k\zz$-graded
complex, and $i\in\zz/k\zz$.
The minimal number $m(D_{i})$
where $D_\stakr$ is a 
$\zz/k\zz$-graded 
complex, chain equivalent to $C_\stakr$,
is denoted by $\mu_{i}(C_\stakr)$.
\end{defi}

\subsection{Lower bounds from local coefficient homology}
\lb{su:znloccoef}

Similarly to section \ref{su:loc-coef}
we have the following estimates for the invariants $\mu_i$
of $\zz/k\zz$-graded complexes.
The proof of the next Proposition is similar 
to \ref{p:bn}.

\bepr\lb{p:z2estimate}
Let $\r:G\to \GL(r, \ff)$ a representation. 
Suppose that there is a group epimorphism $\pi_1(X) \to G$, then 
\begin{equation*}
\mu_{i}(\CCCC^\circ_\stakr(X)) \geq 
\sum_{s\in i(k)} \b_{s}(X,\r). 
\ \qquad \qs
\end{equation*}
\enpr

\subsection
{Invariant $\mu_1$ in the $k$-graded case}
\lb{su:mu1-ngraded}

The previous theorem implies the following lower bounds 
for $\mu_{i}(\CCCC^\circ_\stakr(X))$ in terms of 
the invariants $d(G), \d(G)$.

\beth\lb{t:mu+dg}
Suppose that there exists a group  epimorphism from $\pi_1(X)$ to 
a finite and non-trivial group $G$. Then 
\been\item
$ 
\mu_{1}(\CCCC^\circ_\stakr(X)) \geq \max(\d(G)-1, 1).
$
\item
If $G$ is simple or solvable 
we have 
$\mu_{0}(\CCCC^\circ_\stakr(X))
+
\mu_{1}(\CCCC^\circ_\stakr(X)) \geq d(G).
$
\enen
\enth
\Prf 
We need only to prove that 
$\mu_1(\CCCC^\circ_\stakr(X))\geq 1$.
To this end, observe that if $\mu_1(\CCCC^\circ_\stakr(X))=0$,
then the homology of $X$ in degree 1 with all 
local coefficients vanish, which imply 
$\d(G)=1$, then $G$ is cyclic, $b_1(G)=1$, 
which leads to a contradiction. $\qs$

For the case of infinite groups we have the following
result.
\beth\label{geq1}
Suppose that there exists a group epimorphism from $\pi_1(X)$ to an infinite group $G$.  
Then 
$\mu_{1}(\CCCC^\circ_\stakr(X))\geq 1 $.
\enth
\Prf
Let $D_\stakr$ be a $\zz/k\zz$-graded chain 
complex equivalent to
$\CCCC^\circ_\stakr(X)$.
The module $\wh\L=H_0(\CCCC_\stakr(X))$ is a submodule 
of $H_0(\CCCC^\circ_\stakr(X))$. The condition
$D_1=0$ would imply that $H_0(\CCCC_*^\circ(X))$ and $\wh\L$ are submodules 
of a free $\LL$-module $\CCCC_0^\circ(X)$, and this is 
impossible when $G$ is infinite by Lemma \ref{l:non-embedded}. $\qs$

This theorem holds for $\widehat{\Lambda} \otimes_{\mathbb Z} {\mathbb Q}$ by the same argument.  
This estimate can be improved in the case when 
$k-2\geq \dim X$. Note that in this case 
the sum in the right hand side of 
\rrf{f:svertka} contains only one term for every 
$i$.

\beth\lb{t:big-grade}
Assume that $\dim X\leq k-2$, and there exists a group epimorphism from $\pi_1(X)$ to a finite group $G$.
Then 
$$\mu_1(\CCCC^\circ_*(X))\geq \d(G).$$
\enth
\Prf
Similarly to Section \ref{su:finite-groups}
it suffices to prove that 
$$
\mu_1(\CCCC^\circ_*(X))
\geq
\b_1(X,\r)+1 
$$
for any representation $\r$ such that $H_0(X,\r)=0$.
Let $$D_*=\{\ldots \la D_{-2}\arrl {\pr_{-1}}
D_{-1} \arrl {\pr_{0}} D_0 \la\ldots \}$$
be a $\zz/k\zz$-graded complex, chain equivalent to 
$K_*=\CCCC^\circ_*(X)$.
The chain complex $D_*$ does not necessarily vanish in any degree, and
the argument which we used in the 
proof of the Proposition \ref{p:bn1} can not be applied immediately.

Let 
$D_*\arrr \phi K_*\arrr \psi D_*$
be the mutually inverse chain equivalences.
Since $k\geq \dim X+2$, the
chain complex $K_*$ vanishes in degree  $-1$, 
hence the map $\psi\circ\phi:D_{-1}\to D_{-1}$ is
equal to $0$. The existence of chain homotopy
from $\psi\circ\phi$ to $\Id$ implies that the submodule 
$\Ker\pr_{-1} = \Im \pr_{0}$ is a 
direct summand of $D_{-1}$, hence a projective 
$\wh\L[G]$-module. Let us denote it by $L$.
The $\zz/k\zz$-graded chain complex $D_*$ contains a ($\zz$-graded) subcomplex 
$$
D'_*
=
\{0 \la  L \arrl {\pr_0} 
D_0 \la D_1 \la D_2 \la 0 \},
$$
with
$H_0(D'_*)\approx \wh\L, \ H_{-1} (D'_*)=0,
H_{1} (D'_*)\approx H_{1} (D_*).$
Denote by $S\sbs \nn$ the multiplicative subset 
of all numbers $t$, such that $gcd(t, |G|)=1$.
The module $S^{-1} L$ is free 
by a fundamental result of R. Swan
(see \cite{Roggenkamp}, \S 5).
Thus the chain complex
$D''_*=S^{-1}D'_*$ is free over $S^{-1}\wh\L[G]$; put
$$\a=m(D''_{-1}), \ 
\b= m(D''_{0}), \
\g= m(D''_{1}).
$$

Let $p$ be a prime divisor of $|G|$,
and $\r:G\to \GL(r,\ff_p)$ be a representation.
The homology of the complex 
$K_*\tens{\r}\ff^r_p$ 
is isomorphic to that of 
$D''_*\tens{\r}\ff^r_p$
in degrees $-1,0,1$.
Therefore the argument proving 
Proposition \ref{p:bn1} 
applies here as well, 
and the proof of the Theorem is complete. $\qs$

\section{Estimates for the number of closed orbits}
\label{s:estimates}

We proceed to the estimates of the number of periodic orbits of 
a Hamiltonian isotopy induced by a non-degenerate 
$1$-periodic Hamiltonian $H$ on a closed connected symplectic 
manifold $M$.
We denote by 
$\wi M\to M$ a regular covering with a structure group
$G$. Put
\begin{equation}\lb{f:c-m}
\CCCC_*(M)=
C_*(\wi M)\tens{\zz[G]}
\Lambda^{\mathbb Z[G]}_{(M,\omega)}.
\end{equation}
Denote by $N$ the minimal Chern number of $M$.

\bede\label{d:mui-m}
If $N=0$, put
\begin{equation}\lb{f:mi-m}
 \mu_i(\wi M)=\mu_i(\CCCC_*(M)).
\end{equation}
(see Definition \ref{d:mu}; here $i\in\nn$).

If $N>0$ 
consider the $\zz/2N\zz$-graded chain complex
$C^\circ_*(\wi M)$ (see \rrf{f:svertka}) and put
\begin{equation}\lb{f:mi-mN}
 \mu_i(\wi M)=\mu_i(\CCCC^\circ_*(M)).
\end{equation}
(see Definition \ref{d:mumu}; here $i\in\zz/2N\zz$).

\end{defi}

The numbers $\mu_i(\wi M)$ are obviously homotopy
invariants of $M$ and the chosen covering $\wi M\to M$.

\subsection{ The spherical Calabi-Yau case}
\label{su:CYY}

We consider here symplectic manifolds 
$M$ with $c_1(M)(A)=0$ for every $A\in\pi_2(M)$.
In this case every contractible periodic orbit 
$\g$ has a well-defined index $i(\g)\in\zz$.
The Floer chain complex $\wi{CF}_*$
is a  $\zz$-graded \fbfg~ \cco~
over the ring $\Lambda^{\mathbb Z[G]}_{(M,\omega)}$, generated in degree 
$k$ by contractible 
periodic orbits of the Hamiltonian vector field 
of index $k$, and we have 
$$\wi{CF}_*\sim \CCCC_{*+n}(M),
\ \ {\rm where } \ \ 
\dim M = 2n.$$

Denote by $p_k$ the number of contractible periodic 
orbits of period $k$. The results of the previous sections
imply the following lower bound:
$$
p_{i-n}\geq \mu_i(\wi M).
$$
Applying the results of Section 
\ref{s:z-graded} we obtain the following lower bounds  
for $p_i$:
\bepr\label{p:CY-b-i}
For any field $\ff$ and any representation
$\r:G\to\GL(r,\ff)$ we have
$$
p_{i-n}\geq \frac 1r b_i(M, \r).
$$
\enpr
\noindent
Theorem \ref{t:mu1} and  Corollary
\ref{c:mu}
imply some stronger 
lower bounds for $i=1$.
\beth\label{t:CY-index1}
\been\item
If $\pi_1(M)$ is non-trivial,
then $p_{1-n}\geq 1$.
\item
If $\pi_1(M)$ has an 
epimorphism onto a finite group
$G$, then
\been\item
$p_{1-n}\geq \delta(G)$,
\item
$p_{1-n}\geq d(G)$ if $G$ is solvable or simple,
\item $p_{1-n}\geq 2$ if $G$ is not cyclic.
\enen
\enen 
\enth
Let us proceed to the lower bounds for $p_{2-n}$.
Applying Corollaries \ref{c:mu2_first} and \ref{c:g-finite-mu2}
we obtain the following result.
\beth\label{t:CY-index2}
Assume that $\pi_1(M)$ is finite and the homomorphism 
$\pi_1(M)\to G$ is an isomorphism.
Then
\been\item
$p_{2-n}\geq \delta(\pi_1(M))-b_1(M,\ff)+b_2(M,\ff)$\ \ 
for any field $\ff$.
\item If $G$ is perfect, then $p_{2-n}\geq b_2(M,\ff)+2$.
\item
$p_{2-n}\geq \mu_2(\pi_1(M))$.
\enen
\enth

\bere\label{r:panov}
A recent result of Joel Fine and Dmitry Panov
\cite{FinePanov}
asserts that for every finitely presented group $G$
there exists a symplectic manifold $M$ of dimension $6$
with the fundamental group $G$ and $c_1(M)=0$.
\enre

\pa
\subsection{ The weakly monotone case}
\label{su:w-mono}

Let us denote by $p$ the 
total number of the periodic orbits of the Hamiltonian vector
field.
We have a $\zz/2N\zz$-graded 
chain complex $\wi{CF}_\stakr$, generated by periodic 
orbits, such that
$$
\wi{CF}_\stakr\sim\CCCC^\circ_{\stakr+n}(M),
\ \ {\rm where } \ \ 
2n=\dim M, 
$$
(see Theorem \ref{t:w-monoton-cov}).
Therefore
$$ 
p_{i-n}\geq \mu_{i}(\wi M) \ \text{\rm \ for \ } \ i \in \zz/2N\zz.
$$
Applying the results of the section
\ref{s:z2-graded}
we obtain the following.
\beth\lb{t:monotonBetti}
For every representation 
$\r:G\to \GL(r, \ff)$ we have
\begin{equation*}
p_{i-n} \geq 
\frac 1r \Big(\sum_{s \equiv  i (2N)} b_{s}(M,\r)\Big) \ \text{\rm \ for \ } \ i \in \zz/2N\zz. 
\end{equation*}
\enth
As for the number $p_{1-n}$ we have the following.

\beth\lb{t:monoton}
\been\item 
If $\pi_1(M)$ is non-trivial, then
$p_{1-n}\geq 1$.
\item
If $\pi_1(M)$ has an epimorphism onto a finite group $G$, 
then 
\been\item
$p_{1-n}\geq \max(1, \d(G)-1)$,
and $p\geq \d(G)$.
\item 
if $G$ is simple or solvable,
then $p\geq d(G)$,
\item 
if $G$ is not cyclic, then $p\geq 2$.
\enen
\enen 
\enth

For the manifolds where the minimal Chern number 
$N$ is strictly greater than $n=\dim M/2$ we have 
the following improvement of Theorem \ref{t:monoton}
(the proof  follows from Theorem \ref{t:big-grade}):

\beth\lb{t:monoton-bigChern}
Let $N\geq n+1$.
Assume that $\pi_1(M)$
has an epimorphism onto a finite group $G$.
Then
\been\item
$p_{1-n}\geq \d(G)$.
\item 
If $G$ is simple or solvable,
then $p_{1-n}\geq d(G)$.
\item 
If $G$ is not cyclic, then $p_{1-n}\geq 2$.
\enen
\enth

\subsection{ The general case}
\label{su:gen-case}

Let $M^{2n}$ be an arbitrary  closed connected symplectic manifold.
Theorem \ref{t:gen-case-cover}
together with Proposition \ref{p:betti1}, Theorem \ref{geq1} 
implies the following result.  
As we noted, Theorem \ref{geq1} holds for $\widehat{\Lambda} \otimes_{\mathbb Z} {\mathbb Q}$.   
\beth\label{t:inf-gen-case}
Assume that $\pi_1(M)$ is infinite.
Then $p_{1-n}\geq 1$.
\enth

\section{Acknowledgments}
\label{s:ack}

The authors started this work during
the visit of the second author to the 
Kyoto RIMS in May 2013, and continued 
it during the visit of the second author 
to the Nantes University in November 2013.
It was finalized during the visit of the second 
author to the Kavli IPMU, the University of Tokyo
in April 2014.  

The first author thanks Laboratoire Jean Leray, Universit\'e de Nantes 
for the financial support and its hospitality.      
The second author gratefully acknowledges 
the support of the Kyoto RIMS, the Program for Leading Graduate
Schools, MEXT, Japan, and of the Kavli IPMU, 
the University of Tokyo, and thanks the Kyoto RIMS and 
the Kavli IMPU for hospitality. 
The second author thanks A. Lucchini for sending the papers
\cite{Lucchini90}, \cite{Lucchini96}.

\end{document}